\documentclass[11pt]{article}

\usepackage{color}
 \usepackage[T1]{fontenc}
\usepackage[utf8]{inputenc} 

\newtheorem{myproposition}{Proposition}[section]
\newtheorem{mytheorem}[myproposition]{Theorem}
\newtheorem{mylemma}[myproposition]{Lemma}

\newtheorem{myobservation}[myproposition]{Observation}

\usepackage{amsfonts}
\usepackage{amsmath}
\def\gr{\mathcal{G}}
\def\zet{\mathbb{Z}}
\newcommand{\qed}{\hfill \rule{.1in}{.1in}}

\makeatletter
\def\imod#1{\allowbreak\mkern10mu({\operator@font mod}\,\,#1)}
\makeatother

\begin{document}

\title{Note on group irregularity strength of disconnected graphs}


\author{Marcin Anholcer$^1$, Sylwia Cichacz$^{2}$,  Rafa{\l} Jura$^2$, Antoni Marczyk$^2$\\
$^1$Pozna\'n University of Economics and Business\\
$^2$AGH University of Science and Technology Krak\'ow, Poland}



\maketitle

%



\begin{abstract}
We investigate the \textit{group irregularity strength} ($s_g(G)$) of graphs, i.e. the smallest value of $s$ such that taking any Abelian group $\gr$ of order $s$, there exists a function $f:E(G)\rightarrow \gr$ such that the sums of edge labels at every vertex are distinct. So far it was not known if $s_g(G)$ is bounded for disconnected graphs. In the paper we we present some upper bound for all graphs.  Moreover we give the exact values and bounds on $s_g(G)$ for disconnected graphs without a star as a component. 
\end{abstract}

\section{Introduction}

It  is a well known fact that in any simple graph $G$
there are at least two vertices of the same degree. The situation changes if we consider an edge labeling $f:E(G)\rightarrow \{1,\ldots,s\}$ and calculate weighted degree (or weight) of each vertex $x$ as the sum of labels of all the edges incident to $x$. The labeling $f$ is called \textit{irregular} if the weighted degrees of all the vertices are distinct. The smallest value of $s$ that allows some irregular labeling is called \textit{irregularity strength of $G$} and denoted by $s(G)$.

The problem of finding $s(G)$ was introduced by Chartrand et al. in \cite{ref_ChaJacLehOelRuiSab1} and investigated by numerous authors  \cite{ref_AigTri,ref_AmaTog,ref_FerGouKarPfe,ref_Leh,ref_Tog1}. Best published general result due to Kalkowski et al. (see \cite{ref_KalKarPfe1}) is $s(G)\leq 6n/\delta$. It was recently improved by Przyby{\l}o (\cite{ref_Prz3}) for dense graphs of sufficiently large order ($s(G)\leq (4+o(1))n/\delta+4$ in this case). 

Jones combined the concepts of graceful labeling and modular edge coloring into labeling called  a \emph{modular edge-graceful labeling} (\cite{ref_Jon,ref_JonKolOkaZha2, ref_JonZha}). He defined the \textit{modular edge-gracefulness of graphs} as the smallest integer $k(G)=k \geq n$ for which there exists an edge labeling $f:E(G)\rightarrow \zet_k$ such that the induced vertex labeling $f^\prime : V (G)\rightarrow \zet_k$ defined by
$$
f^\prime(u) =\sum_{v\in N(u)}f(uv)\mod k
$$
is one-to-one.

Assume $\gr$ is an Abelian group  of order $m\geq n$ with the operation denoted by $+$ and identity element $0$. For convenience we will write $ka$ to denote $a+a+\ldots+a$ (where element $a$ appears $k$ times), $-a$ to denote the inverse of $a$ and we will use $a-b$ instead of $a+(-b)$. Moreover, the notation $\sum_{a\in S}{a}$ will be used as a short form for $a_1+a_2+a_3+\ldots$, where $a_1, a_2, a_3, \ldots$ are all the elements of the set $S$. Recall that any group element $\iota\in\gr$ of order 2 (i.e., $\iota\neq 0$ such that $2\iota=0$) is called \emph{involution}. 

The order of an element $a\neq 0$ is the smallest $r$ such that $ra=0$. It is well-known by Lagrange Theorem that $r$ divides $|\gr|$ \cite{ref_Gal}. Therefore every  group of odd order has no involution.

We consider edge labeling $f:E(G)\rightarrow \gr$ leading us to the weighted degrees defined as the sums (in $\gr$):
$$
w(v)=\sum_{v\in N(u)}f(uv)
$$

The concept of $\gr$-irregular labeling is a generalization of modular edge-graceful labeling. In both cases the labeling $f$ is called \textit{$\gr$-irregular} if all the weighted degrees are distinct. However, the \textit{group irregularity strength} of $G$, denoted $s_g(G)$, is the smallest integer $s$ such that for every Abelian group $\gr$ of order $s$ there exists $\gr$-irregular labeling $f$ of $G$. Thus the following observation is true.

\begin{myobservation}[\cite{ref_AnhCic2}]
For every graph $G$ with no component of order less than $3$, $k(G)\leq s_g(G)$.
\end{myobservation}

The following theorem, determining the value of $s_g(G)$ for every connected graph $G$ of order $n\geq 3$, was proved by Anholcer, Cichacz and Milani\v{c} \cite{ref_AnhCic1}.

\begin{mytheorem}[\cite{ref_AnhCic1}]\label{AnhCic1}
Let $G$ be an arbitrary connected graph of order $n\geq 3$. Then
$$
s_g(G)=\begin{cases}
n+2,&\text{if   } G\cong K_{1,3^{2q+1}-2} \text{  for some integer   }q\geq 1,\\
n+1,&\text{if   } n\equiv 2 \imod 4 \wedge G\not\cong K_{1,3^{2q+1}-2} \text{  for any integer   }q\geq 1,\\
n,&\text{otherwise.}
\end{cases}
$$
\end{mytheorem}

In \cite{ref_Jon} it was proved in turn that for every connected graph $G$ of order $n\geq 3$ 
$$
k(G)=\left\{\begin{array}{lll}
n,&\text{if}& n \not \equiv 2 \imod 4, \\
n+1,&\text{if}& n\equiv 2 \imod 4.
\end{array}
\right.
$$

In order to distinguish $n$ vertices in arbitrary (not necessarily connected) graph we need at least $n$ distinct
elements of $\gr$. However, $n$ elements are not always enough, as shows the
following lemma.

\begin{mylemma}[\cite{ref_AnhCic2}]\label{lemma_below}
Let $G$ be a graph of order $n$. If $n \equiv 2 \imod 4$, then there is no $\gr$-irregular labeling of $G$ for any Abelian group $\gr$ of order $n$.
\end{mylemma}

Anholcer and Cichacz considered the group irregularity strength of disconnected graphs in \cite{ref_AnhCic2}.
\begin{mytheorem}[\cite{ref_AnhCic2}]\label{inwolucja}
Let $G$ be a graph of order $n$ with no component of order less than $3$ and with all the bipartite components having both color classes of even order. Let $s=n+1$ if $n\equiv 2 \imod 4$ and $s=n$ otherwise. Then:
$$
\begin{array}{lll}
s_g(G)=n,& \text{if} &n \equiv 1\imod 2,\\
s_g(G)=n+1,& \text{if} &n\equiv 2 \imod 4,\\
s_g(G)\leq n+1, & \text{if} &n\equiv 0 \imod 4.
\end{array}
$$
 Moreover for every integer $t\geq s$ there exists a $\gr$-irregular labeling of $G$ for every Abelian group $\gr$ of order $t$ with at most one involution $\iota$. 
\end{mytheorem}

\begin{mytheorem}[\cite{ref_AnhCic2}]\label{modulo}
Let $G$ be a graph of order $n$ having neither component of order less than $3$ nor a $K_{1,2u+1}$ component for any integer $u\geq 1$. Then:

$$
\begin{array}{lll}
k(G)=n,& \text{if} &n \equiv 1\imod 2,\\
k(G)=n+1,& \text{if} &n\equiv 2 \imod 4,\\
k(G)\leq n+1, & \text{if} &n\equiv 0 \imod 4.
\end{array}
$$
Moreover for every odd integer $t\geq k(G)$ there exists a $\zet_t$-irregular labeling of $G$.
\end{mytheorem}

In this paper we give an upper bound for group irregularity strength of all graphs.
Moreover we give the exact values and bounds on $s_g(G)$ for disconnected graphs with no star components.

\section{Main results}

The first natural question is whether the group irregularity strength is finite for arbitrary graph with no components of order less than $3$.

\begin{mytheorem}\label{main2}
Let $G$ be a graph of order $n$ having $m$ components, none of which has order less than $3$ and let $p$ be the smallest prime number greater than $2^{n-m-1}$. Then $s_g(G)\leq p$.
\end{mytheorem}
\noindent\textbf{Proof.}
Note that $n\leq p$ and $p$ is odd. Since $p$ is prime, there exists only one (up to isomorphism) Abelian group $\gr$ of order $p$, namely $\gr\cong \zet_p$. If all the components of $G$ have order $3$ then we are done by Theorem~\ref{modulo}. Therefore $G$ has at least one component $H$ such that $|V(H)|\geq4$.  Let $F$ be a spanning forest of $G$. Thus $F$ has $n-m$ edges $e_0,e_1,\ldots,e_{n-m-1}$. Assume without loss of generality that $e_0\in E(H)$. Let $f\colon E(G) \rightarrow \zet_p$ be defined as follows 
$$
\begin{array}{lll}
f(e_i)=2^{i-1},& \text{for} &i=1,2,\dots n-m-1,\\
f(e)=0,& \text{for} &e\in \{e_0\}\cup E(G)\setminus E(F).
\end{array}
$$
We have $\sum_{i=0}^{k-1}2^i=2^k-1<2^k$ for any integer $k$. Therefore the maximum weighted degree is smaller than $2^{n-m-1}$. Moreover the unique (additive) decomposition of any natural number into powers of $2$ implies that $f$ is $\gr$-irregular.~\qed

Given any two vertices $x_1$ and $x_2$ belonging to the same connected component of $G$, there exist walks from $x_1$ to $x_2$. Some of them may consist of even number of vertices (some of them being repetitions). We are going to call them \textit{even walks}. The walks with odd number of vertices will be called \textit{odd walks}. We will always choose the shortest even or the shortest odd walk from $x_1$ to $x_2$.

We start with $0$ on all the edges of $G$. Then, in every step we will choose $x_1$ and $x_2$ and add some labels to all the edges of chosen walk from $x_1$ to $x_2$. To be more specific, we will add some element $a$ of the group to the labels of all the edges having odd position on the walk (starting from $x_1$) and $-a$ to the labels of all the edges having even position. It is possible that some labels will be modified more than once, as the walk does not need to be a path. We will denote such situation with $\phi_e(x_1,x_2)=a$ if we label the shortest even walk and $\phi_o(x_1,x_2)=a$ if we label the shortest odd walk. Observe that putting $\phi_e(x_1,x_2)=a$ results in adding $a$ to the weighted degrees of both $x_1$ and $x_2$, while $\phi_o(x_1,x_2)=a$ means adding $a$ to the weighted degree of $x_1$ and $-a$ to the weighted degree of $x_2$. In both cases the operation does not change the weighted degree of any other vertex of the walk. Note that if some component $G_1$ of $G$ is not bipartite, then for any vertices $x_1,x_2\in G_1$ there exist both even and odd walks.

We are going to use the following theorem, proved in \cite{ref_Zeng}.

\begin{mytheorem}[\cite{ref_Zeng}]\label{thmZSP}
Let $s=r_1+r_2+\ldots+r_q$ be a partition of the positive integer $s$, where $r_i\geq 2$ for $i=1,2,\ldots,q$. Let $\gr$ be an Abelian group of order $s+1$ .
Then the set $\gr\setminus\{0\}$ can be partitioned into pairwise disjoint subsets $A_1, A_2, \ldots, A_q$ such that for every $1\leq i\leq q$, $|A_i|=r_i$ with
$\sum_{a\in A_i}{a}=0$ if and only if $|\gr|$ is odd or $\gr$ contains exactly three involutions.
\end{mytheorem}

From the above Theorem~\ref{thmZSP} we easily obtain the following observation:
\begin{myobservation}\label{oZSP}
Let $s=r_1+r_2+\ldots+r_q$ be a partition of the positive odd integer $s$, where $r_i\geq 2$ for $i=2,3,\ldots,q$. Let $\gr$ be an Abelian group of order $s$.
Then the set $\gr$ can be partitioned into pairwise disjoint subsets $A_1, A_2, \ldots, A_q$ such that for every $1\leq i\leq q$, $|A_i|=r_i$ with
$\sum_{a\in A_i}{a}=0$.~\qed
\end{myobservation}

Using the simmilar method as in the proof of Theorem~\ref{inwolucja}, we can obtain the following lemma.
\begin{mylemma}\label{lemZSP}
Let $G$ be a graph of order $n$ having no $K_{1,u}$ components for any integer $u\geq 0$. Then for every odd integer $t\geq n$ and for every Abelian group $\gr$ such that $|\gr|=t$, there exists a $\gr$-irregular labeling.
\end{mylemma}  

\noindent\textbf{Proof.}

We are going to divide the vertices of $G$ into triples and pairs. Let $p_1$ be the number of bipartite components of $G$ with both color classes odd, $p_2$ with both classes even and $p_3$ with one class odd and one even. Let $p_4$ be the number of remaining components of odd order and $p_5$ - the number of remaining components of even order. The number of triples equals to $2p_1+p_3+p_4$. The remaining vertices form the pairs.

By Observation~\ref{oZSP}, the elements of $\gr$ can be partitioned into $2l+1$ triples $B_1,B_2,\ldots,B_{2l+1}$ and $m$ pairs $C_1,C_2,\ldots,C_m$, where $l=\lfloor(2p_1+p_3+p_4)/2\rfloor$ and $m=(t-6l-3)/2$, such that $\sum_{x\in B_i}x=0$ for $i=1,\dots,2l+1$ and $\sum_{x\in C_j}x=0$ for $j=1,\dots,m$. Observe that $l\geq 0$, $m\geq 0$ and $2l+1\geq 2p_1+p_3+p_4$.

Let $B_i=\{a_i,b_i,c_i\}$ for $i=1,2,\ldots,2l+1$ and let $C_j=\{d_j,-d_j\}$ for $j=1,\ldots,m$. It is easy to observe that for a given element $g\in\gr$ not belonging to any triple, we have $(g,-g)=C_j$ for some $j$. 


Let us start the labeling. For both vertices and labels, we are numbering the pairs and triples consecutively, in the same order as they appear in the labeling algorithm described below, every time using the lowest index that has not been used so far (independently for the lists of couples and triples).

Given any bipartite component $G$ with both color classes even, we divide the vertices of every color class into pairs $(x_j^1,x_j^2)$, putting
$$
\phi_o(x_j^1,x_j^2)=d_j
$$
for every such pair. We proceed in similar way in the case of all the non-bipartite components of even order, coupling the vertices of every such component in any way.

If both color classes of a bipartite component are of odd order, then they both have at least $3$ vertices. We choose three of them, denoted with $x_j$, $y_j$ and $z_{j}$, in one class and another three, $x_{j+1}$, $y_{j+1}$ and $z_{j+1}$, in another one and we put
$$
\begin{array}{l}
\phi_e(x_j,z_{j+1})=a_j,\\
\phi_e(y_j,z_{j+1})=b_j,\\
\phi_e(z_j,z_{j+1})=c_j,\\
\phi_e(x_{j+1},z_{j})=a_{j+1},\\
\phi_e(y_{j+1},z_{j})=b_{j+1},\\
\phi_e(z_{j+1},z_{j})=c_{j+1}.\\
\end{array}
$$
We proceed with the remaining vertices of these components as in the case when both color classes are even.

In the case of  non-bipartite components of odd order we choose three vertices. We put
$$
\begin{array}{l}
\phi_e(x_j,z_j)=a_j,\\
\phi_e(y_j,z_j)=b_j,\\
\phi_e(z_j,z_j)=c_j.
\end{array}
$$

Finally for bipartite components of odd order we choose four vertices $x_j$, $y_j$, $z_j$ and $v$ ($v$ belongs to the even color class and three other vertices to the odd one). We put
$$\begin{array}{l}
\phi_e(x_j,v)=a_j,\\
\phi_e(y_j,v)=b_j,\\
\phi_e(z_j,v)=c_j.
\end{array}
$$

The labeling defined above is $\gr$-irregular. Indeed, in the $j^{th}$ triple of vertices the weights are equal to $w(x_j)=a_j$, $w(y_j)=b_j$ and $w(z_j)=c_j$ and in the $j^{th}$ pair we have $w(x_j^1)=d_j$ and $w(x_j^2)=-d_j$. Eventually, at least one of the triples of labels remains unused.
\qed

The following theorem easily follows from the above Lemmas~\ref{lemma_below} and~\ref{lemZSP}.
\begin{mytheorem}\label{main}
Let $G$ be a graph of order $n$ having no $K_{1,u}$ components for any integer $u\geq 0$. Then:
$$
\begin{array}{lll}
s_g(G)=n,& \text{if} &n \equiv 1\imod 2,\\
s_g(G)=n+1,& \text{if} &n\equiv 2 \imod 4,\\
s_g(G)\leq n+1, & \text{if} &n\equiv 0 \imod 4.
\end{array}
$$\qed
\end{mytheorem}

We will consider now some families of disconnected graphs of order $n\equiv 0 \imod 4$ for which 
$s_g(G)=n$. 

\begin{myproposition}
Let $G$ be a graph of order $n\equiv 4\imod 8$ with no component of order less than $3$ and with all the bipartite components having both color classes of even order. Then $s_g(G)=n$.
\end{myproposition}
\noindent\textbf{Proof.}

Let $\gr$ be an Abelian group of order $n$. Since the order of $\gr$ is even there is at least one involution in $\gr$. If there is exactly one involution, then we are done by Theorem~\ref{inwolucja}.  Thus we can assume that $\gr$ has more than one involution.  Observe that $n=2^2(2\alpha+1)$ for some integer $\alpha$, therefore by fundamental theorem of finite Abelian groups we obtain that $\gr$ has exactly three involutions $\iota_1, \iota_2,\iota_3$.

Let $p_1$ be the number of components of odd order, $p_2$  be the number of components of even order. 

Assume first $p_2>0$. Then there exists a component $H$ of even order $|H|\geq 4$. Note that there exist vertices $u,v,x,y\in V(H)$ such that there is an odd walk from $u$ to $x$, an even walk from $u$ to $v$ and an even walk from $u$ to $y$ (if $H$ is bipartite, we take $u$ and $x$ from one color class and $v$ and $y$ from another, what is always possible, since in this case $H$ both color classes have even order). By Theorem~\ref{thmZSP}, the set of the elements of $\gr\setminus\{0\}$ has partition into $p_1+1$ triples $B_1,B_2,\ldots,B_{p_1+1}$ and $m$ pairs $C_1,C_2,\ldots,C_m$ where $m=(n-3p_1-4)/2\geq 0$ such that $\sum_{x\in B_i}x=0$ for $i=1,\dots,p_1+1$ and $\sum_{x\in C_j}x=0$ for $j=1,\dots,m$. Let $B_{p_1+1}=\{a_{p_1+1},b_{p_1+1},c_{p_1+1}\}$, without loss of generality we can assume that $a_{p_1+1}=\iota_1$. Put
$$\begin{array}{l}
\phi_o(u,x)=a_{p_1+1},\\
\phi_e(u,v)=b_{p_1+1},\\
\phi_e(u,y)=c_{p_1+1}.
\end{array}
$$
Note that we obtain now $w(u)=0$, $w(v)=-\iota_1=\iota_1$, $w(x)=b_{p_1+1}$ and $w(y)=c_{p_1+1}$. We proceed with the remaining vertices in the same way as in the proof of Lemma~\ref{lemZSP} (we divide $V(G)\setminus \{x,y,u,v\}$ into triples and pairs). 

If $p_2=0$ then by Theorem~\ref{thmZSP}, the set of the elements of $\gr\setminus\{0\}$ has partition into  triples $B_1,B_2,\ldots,B_{p_1-1}$ and $m$ pairs $C_1,C_2,\ldots,C_m$ where $m=(n-3p_1+2)/2\geq 0$ such that $\sum_{x\in B_i}x=0$ for $i=1,\dots,p_1-1$ and $\sum_{x\in C_j}x=0$ for any $j=1,\dots,m$. We set $B_{p_1}=C_m\cup\{0\}$  and proceed in the same way as in the proof of Lemma~\ref{lemZSP}.~\qed\\

For $n\equiv 0 \imod 8$ we have the following result, unfortunately with a stronger assumption on non-bipartite components:
\begin{mytheorem}\label{czworki}
Let $G$ be a disconnected graph of order $n$ with all  components of order divisible by $4$ and all the bipartite components having both color classes of even order. Then $s_g(G)=n$.
\end{mytheorem}
\noindent\textbf{Proof.}
Let $\gr$ be an Abelian group of order $n$. Note that $n \equiv 0 \imod 4$. Since the order of $\gr$ is even there is at least one involution in $\gr$, thus by Theorem~\ref{inwolucja} we can assume that $\gr$ has the set of involutions $I^\star=\{\iota_1,\iota_2,\ldots,\iota_{2^p-1}\}$ for some $p\geq2$.

Obviously $I = I^*\cup\{0\}$ is a subgroup of $\gr$. Note that $\Gamma=\{0,\iota_1,\iota_2,\iota_1+\iota_2\}$ is a subgroup of $I$ as well as a subgroup of $\gr$. If $p=2$ then we define $B_1=\Gamma$. If $p\geq 3$, then there exists a coset decomposition of $I$ into $a_1+\Gamma,a_2+\Gamma,\ldots,a_{2^{p-3}}+\Gamma$ for  $a_j\in I$, $j=1,2,\ldots,2^{p-3}$. Set $B_j=a_j+\Gamma$ for $j=1,2,\ldots,2^{p-3}$. Obviously $\sum_{b\in B_j}b=0$, and moreover for any $b\in B_j$ we have $-b=b$ for $j=1,2,\ldots,2^{p-3}$.

Note that the remaining elements of $\gr$ i.e. the elements of $\gr\setminus I$ can be divided into quadruples of four distinct elements $B_j=\{g^1_j,-g^1_j,g^2_j,-g^2_j\}$ for $j=2^{p-3}+1,2^{p-3}+2,\ldots,|\gr|/4$, none of which being an involution.\\

Let $B_j=\{b^1_j,b^2_j,b^3_j,b^4_j\}$   and  $b_j^3\notin\{b_j^1,-b_j^1\}$ for $j=1,2,\ldots,|\gr|/4$. Let us start the labeling. Given any bipartite component $G$ with both color classes even, we divide the vertices of the component into quadruples $(x_j^1,x_j^2,x_j^3,x_j^4)$ such that, $x_j^1,x_j^2$ are in the same color class, and $x_j^3,x_j^4$ are in the same color class (possibly the same as $x_j^1,x_j^2$ but not necessarily).  We proceed in similar way in the case of all the non-bipartite components.  We are numbering the quadruplets consecutively, starting with $1$.

If there is an involution in $B_j$ then set
$$\begin{array}{l}
\phi_o(x_j^1,x_j^2)=b_j^1,\\
\phi_e(x_j^1,x_j^3)=b_j^2,\\
\phi_e(x_j^1,x_j^4)=b_j^3.
\end{array}
$$
Observe that in that case  $w(x_j^2)=b_j^1+b_j^2+b_j^3=-b_j^4=b_j^4$, $w(x_j^2)=-b_j^1=b_j^1$, $w(x_j^3)=b_j^2$ and $w(x_j^4)=b_j^3$. If there is no involution in $B_j$ then let
$$\begin{array}{l}
\phi_o(x_j^1,x_j^2)=b_j^1,\\
\phi_o(x_j^3,x_j^4)=b_j^3.
\end{array}
$$
Note that we obtain now $w(x_j^1)=-w(x_j^2)=b_j^1$, $w(x_j^3)=-w(x_j^4)=b_j^3$.~\qed\\

The \textit{lexicographic product} or \textit{graph composition} $G \circ H$ of graphs $G$ and $H$ is a graph such that  the vertex set of $G \circ H$ is the Cartesian product $V(G) \times V(H)$ and any two vertices $(u,v)$ and $(x,y)$ are adjacent in $G \circ H$ if and only if either $u$ is adjacent with $x$ in $G$ or $u = x$ and $v$ is adjacent with $y$ in $H$. Note that $G\circ H$ and $H\circ G$ are not isomorphic in general. One can imagine obtaining $G\circ H$ by blowing up each vertex of $G$ into a copy of $H$. For instance $lK_{2r,2r}\cong lK_2\circ \overline{K}_{2r}$. 
		
One can easily see that if $H$ has no isolated vertices and $F$ is a graph of order divisible by $4$, then $s_g(H\circ F)=|H|\cdot|F|$ by the above Theorem~\ref{czworki}. Observe also that if  $H$ has all components of even order then for $G=H\circ\overline{K}_{2r}$ we have $s_g(G)=2r|H|$ for any $r\geq1$. One could ask if we need the assumption on the order of components of $H$. Before we proceed we will need  the following result:
\begin{mytheorem}[\cite{ref_Cic}]\label{mZSP}
Let $s=qr$, where $r\geq 3$ and $\gr$ be an Abelian group of order $s$ such that the number of involutions in $\gr$ is not one.
Then the set $\gr$ can be partitioned into pairwise disjoint subsets $A_1, A_2, \ldots, A_q$ such that for every $1\leq i\leq q$, $|A_i|=r$ with
$\sum_{a\in A_i}{a}=0$.
\end{mytheorem}

\begin{myobservation}\label{last}
Let $H$ be a graph  of order $n$ with no isolated vertices.
If $G\cong H\circ \overline{K}_{2r}$ for some positive integer $r\geq2$ , then
$s_g(G)=2rn$ for $rn$ even and $s_g(G)=2rn+1$ otherwise.
\end{myobservation}
\noindent\textbf{Proof.}

 Obviously $G$ is a graph of order $2nr$ with no component of order less than $3$ and with all the bipartite components having both color classes of even order.  If $nr$ is odd, then $2nr \equiv 2 \imod 4$, hence $s_g(G)=2rn+1$ by Theorem~\ref{inwolucja}.
Therefore we can assume that $2nr \equiv 0 \imod 4$. Let $\gr$ be an Abelian group of order $2nr$.  Since the order of $\gr$ is even there is at least one involution in $\gr$, therefore we can assume that $\gr$ has more than one involution by Theorem~\ref{inwolucja}. The set of the elements of $\gr$ has ampartition into  sets $A_1,A_2,\ldots,A_{n}$ of  order $2r$ such that $\sum_{x\in A_i}x=0$  by  Theorem~\ref{mZSP}.

Let $A_i=\{a^1_i,a^2_i,\ldots,a^{2r}_i\}$ for $i=1,2,\ldots,n$.  Denote the vertices of $G$ corresponding to a vertex $x_i\in V(H)$  by $x_i^1,x_i^2,\ldots,x_i^{2r}$. Let $y\in N_{H}(x_i)$, then $y^1 \in N(x_i^j)$  for $j=1,2,\ldots,2r$. Set $\phi_e(x_i^j,y^1)=a_j^i$ for $j=1,2,\ldots,2r$. One can check that  the weighted degrees of all the vertices are distinct.~\qed

Using the same method as in the proof of Observation~\ref{last} we have the following result.
\begin{myobservation}
Let $H$ be a graph of order $n$ with no isolated vertices and with all the bipartite components having both color classes of even order. If $G\cong H\circ \overline{K}_{2r+1}$ for some positive integer $r\geq2$ , then
$s_g(G)=(2r+1)n$ for $n\not \equiv 2 \imod 4$ and  $s_g(G)=(2r+1)n+1$ otherwise.~\qed
\end{myobservation}

\nocite{*}
\bibliographystyle{amsplain}

\end{document}